\documentclass[11pt]{article}
\usepackage{amsfonts}
\usepackage{amsthm}
\usepackage{multirow}
\usepackage{amsmath}
\usepackage{amssymb}
\usepackage{graphicx}
\usepackage{color}

\usepackage{marginnote}
\usepackage{amstext,amsthm,amssymb,amsmath}
\usepackage{graphicx}
\usepackage{subfigure}

\hfuzz=\maxdimen
\graphicspath{ {./pic/}}
\title{Multicontinuum homogenization. General theory and applications.}
\author{E. Chung\footnote{Department of Mathematics, The Chinese University of Hong Kong, Sha Tin, Hong Kong}, Y. Efendiev\footnote{Department of Mathematics, Texas A\&M University, College Station, TX 77843, USA}, J. Galvis\footnote{Departamento de Matem\'aticas, Universidad Nacional de Colombia, Carrera 45 No. 26-85, Edificio Uriel Guti\'errez, Bogot\'a D.C., Colombia}, W.T. Leung\footnote{Department of Mathematics, City University of Hong Kong, Hong Kong}}

\begin{document}
\maketitle
%==============================

\begin{abstract}

In this paper, we discuss a general framework for multicontinuum 
homogenization. Multicontinuum models are widely used in many
applications and some derivations for these models are established.
In these models, several macroscopic variables at each macroscale point
are defined and the resulting multicontinuum equations are formulated.
In this paper, we propose a general formulation and associated
ingredients that allow performing multicontinuum homogenization.
Our derivation consists of several main parts. In the first part,
we propose a general expansion, where the solution is expressed via
the product of multiple macro variables and associated cell problems. 
The second part consists
of formulating the cell problems. The cell problems are formulated
as saddle point problems with constraints for each continua. Defining
the continua via test functions, we set the constraints as
an integral representation. Finally, substituting the expansion to
the original system, we obtain multicontinuum systems. 
We present an application to the mixed formulation of elliptic equations.
This is a challenging system as the system does not have symmetry.
We discuss the local problems and various macroscale representations
for the solution and its gradient. Using various order approximations,
one can obtain different systems of equations. We discuss 
the applicability of 
multicontinuum homogenization and relate this to high contrast 
in the cell problem. Numerical results are presented. 

\end{abstract}

\section{Introduction}

Many problems have multiscale nature. For example, the flow in porous media occurs in multiscale
media with heterogeneities at multiple scales and high contrast. The simulations of these problems are
often performed on a coarse computational grid, where the grid size is much larger compared to
the scales of heterogeneities. In these simulations, we distinguish two cases in the paper.
The first is the case with no-scale separation and the second is the case with scale separation.
In the first case, approaches use the information within
the entire computational grid or beyond to
derive macroscopic equations. We will not discuss this case in the paper.
 In the second case, the representative volume-based information (which is much smaller compared to the target coarse block)
is used in deriving macroscopic equations.

For the case of no-scale separation, many approaches are developed to account for subgrid effects.
These approaches, e.g., \cite{chung2016adaptiveJCP,GMsFEM13,eh09,hw97,jennylt03,chung2018constraint,oz06_1,altmann2021numerical,
fish2013practical,new2023,contreras2023exponential,abreu2020conservation,abreu2019convergence},
include the construction of multiscale basis functions that are supported in domain larger 
than the target coarse block. Among these approaches, the CEM-GMsFEM \cite{chung2018constraint} is
related to the approaches presented in this paper. In these approaches, the multiscale basis 
functions are computed in oversampled regions. There are several basis functions in each
target coarse block representing different continua effects. These concepts will further 
be used
in multicontinua homogenization.

In the case of scale separation, one uses information in representative volume (which is much smaller compared to the coarse block) to derive effective
properties. The well-known approach includes the homogenization technique
\cite{bensoussan2011asymptotic,blanc2023homogenization,bakhvalov2012homogenisation,allaire1992homogenization}, which is widely used
in many applications. The main idea of this approach is to assume that the solution in each
macroscopic point, can be represented by its average. The homogenization 
method provides a systematic
expansion, which allows for deriving the equations. In this derivation, the small-scale $\epsilon$
is the RVE size. In the derivation, all terms depending on different powers of $\epsilon$ are
separated. The latter is one of the limitations in extending these methods to problems where
the media properties can depend on $\epsilon$ (high-contrast case). 

In this paper, we introduce a general homogenization method, where we assume that the media properties
can have high contrast. In our expansion, we consider that each macroscopic point has several
macroscopic variables associated with it. The macroscopic variables are defined via auxiliary functions
and assumed to be smooth functions. The expansion of the solution via macroscopic variables uses
the solution of local microscopic problems posed in RVE, called solutions of cell problems. These local problems
account for the micro-scale behavior of the solution given certain constraints. These constraints are related
to the definition of macroscopic variables. In particular, 
our first cell problem imposes constraints to represent the 
constants in the average behavior of each continua. 
The consequent cell problems impose
constraints to represent the 
high-order polynomials in the average behavior of each continua.

The multi-continuum homogenization expansion is substituted into the 
fine-scale equations. 
Our next assumptions
include the fact that the integrals in the macroscopic variational formulation can be written in terms
of the integrals over RVE and macroscopic variables are smooth. Using these assumptions, we derive
a system of equations on a coarse grid. The resulting system of equation include additional terms
and can involve higher-order derivatives. These equations share similarities to other models derived
earlier and some terms can be negligible due to high contrast in the media
properties.

Our approaches share some common ingredients with mixture theories 
\cite{rajagopal1995mechanics,truesdell1984thermodynamics,malek20}. In mixture theories,
the conservation of mass and momentum are written for each component. This model can be used
in deriving a general set of macroscopic equations. However, these models do not make any 
specific assumptions on exchange terms. Our models generalize some earlier derived 
model
equations related to works \cite{efendiev2023multicontinuum}, 
dual-permeability models \cite{rubin48,barenblatt1960basic,showalter1991micro,aifantis1979continuum,iecsan1997theory,bunoiu2019upscaling,arbogast1990derivation,bedford1972multi,chai2018efficient,alotaibi2022generalized}
and we establish a general tool for deriving multicontinuum homogenization models.

One of the challenging aspects of multicontinuum homogenization is in formulating cell
problems correctly. We consider large oversampled regions, where we can impose higher-order
polynomial constraints. By imposing averages in each RVE within the oversampled region, our main
test is to guarantee that the solution of the cell problem converges to zero. To achieve
this, one needs a careful formulation of cell problems. For example, in a carefully studied
example of mixed Darcy equations, we show how one can achieve this.
We obtain a generalized Darcy approximation on the coarse grid.

We present numerical examples. In this numerical example, we consider
a mixed formulation between velocity and pressure in Darcy's equation.
Because pressure and velocity are treated separately, their relation at 
the microscale will not necessarily preserve at the macroscale as in the standard
homogenization. We note that there is a
 linear relation between the velocity and the gradient
of the pressure via the multiscale permeability field. 
Because the mixed formulation
is not  symmetric, this causes further challenges
that are addressed in numerical examples when imposing local constraint problems.
Our numerical results show a good convergence as we decrease the mesh size.

The paper is organized as follows. In the next section, we present preliminaries and a simple derivation of multicontinuum homogenization for zero-order equations. In Section 3, we present a general theory for multicontinuum homogenization and also discuss the relation to mixture theory. Section 4 is devoted to mixed-order systems. In Section 5, we present
numerical experiments.

\section{Preliminaries and zero order equation}

To present preliminaries, we consider zero-order equations
(following \cite{blanc2023homogenization}).
We consider the following zero-order equation
\begin{equation}
\label{eq:zero}
A(x) u(x) = f(x),
\end{equation}
where $A(x)$ is a scalar function with multiple scales and high contrast.
For example, we assume $A$ is a periodic function where the period 
consists of two distinct regions
with highly varying coefficients. We denote by $\psi_i$ the characteristic function for the region
$i$, called the $i$th continua.

%Standard homogenization is the following. 
%\[
%u(x)={f(x)\over A(x)}.
%\]
%Averaging over each RVE, we get 
%\[
%U(x)=f(x) \langle {1\over A(x)}\rangle_{R_\omega}.
%\]
%Consequently,
%\[
%\langle {1\over A(x)} \rangle_{R_\omega}^{-1} U(x) = f(x).
%\]

\begin{figure}
\centering
\includegraphics[scale=0.35]{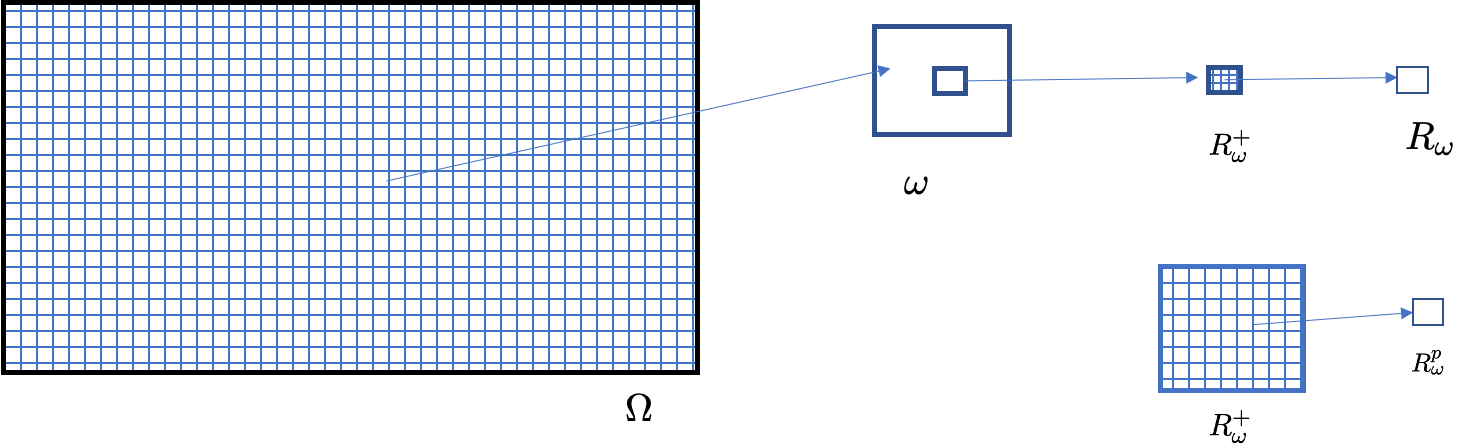}
\caption{Illustration}
\label{fig:ill}
\end{figure}

It is assumed
that the problem is solved on a computational grid consisting
of grid blocks, denoted $\omega$, that are much larger than
heterogeneities.
We assume some type of periodicity within each
computational block represented by Representative 
Volume Element $R_\omega$ that corresponds to 
a computational element $\omega$ (see Figure \ref{fig:ill})
(more precise meaning will be defined later).
We assume that within each $R_\omega$, 
there are several distinct average states
 (known as multicontinua).
We denote the characteristic function for the
continuum $i$ within $R_\omega$ by $\psi_i^\omega$
($\omega$ will be omitted since local computations are restricted to a
coarse block), 
i.e., $\psi_i=1$ within
continuum $i$ (can be irregularly shaped regions consisting of several
parts, in general) and $0$ otherwise.
We introduce oversampled $R_\omega^+$ that contains several $R_\omega^p$'s,
where $p$ denotes different $R_\omega$'s.
We denote the central (target) RVE by, simply, $R_\omega$.
We denote $\psi^p_i$, the characteristic function for 
$R_\omega^p$ and will omit the index $p$ for simplicity if 
it is clear which region we are referring to.

We consider the expansion of $u$ in each RVE as
(for simplicity, we use equal sign instead of approximation)
\begin{equation}
\label{eq:u11}
u=\phi_i U_i,
\end{equation}
where $\phi_i$ is a microscopic function for each $i$ and $U_i(x)$ is
a smooth function for each $i$.
The summation over repeated indices is taken. To obtain the microscopic function $\phi_i$,
we formulate the following cell problems in each RVE within $\omega$ and use $y$ dependence to denote
microscopic nature:
\begin{equation}
\begin{split}
A(y) \phi_i(y) = D_{ij}\psi_j\ \ \text{in} \ R_\omega\\
 \int_{R_\omega} \phi_i \psi_j  = \delta_{ij} \int_{R_\omega} \psi_j,
\end{split}
\end{equation}
where $D_{ij}$ are constants and can be shown that $D_{ij}\psi_j=C_i \psi_i$. Moreover, it can be easily
computed that 
\[
C_j ={\int_{R_\omega} \psi_j  \over \int_{R_\omega} \psi^2_j A^{-1} }.
\]

Next, we derive macroscopic equations. For this, we 
first write an integral form (for any test function $v$)
\begin{equation}
\begin{split}
\int_D f v &= \int_D A(x) u(x) v(x)= \sum_\omega \int_\omega A(x) u(x) v(x) \\
&\approx \sum_\omega {|\omega|\over |R_\omega|} \int_{R_\omega} A(y) u(y) v(y).
\end{split}
\end{equation}
Substituting $u$ from (\ref{eq:u11}) 
into the equation and writing $v=\phi_i V_i$, we get
\begin{equation}
\begin{split}
\int_{R_\omega} A(y) u(y) v(x)dy \approx U_i(x_\omega) V_j(x_\omega) \int_{R_\omega} A(y)\phi_i(y) \phi_j(y) \; dy,
\end{split}
\end{equation}
where $x_\omega$ is a mid-point of $R_\omega$. We will omit
the microscopic dependence of macroscale variables (e.g., $U_i$) and
simply use $U_i$ notation
\begin{equation}
\begin{split}
\int_{R_\omega} A(y) u(y) v(x)dy =U_i V_j\int_{R_\omega} A(y)\phi_i(y) \phi_j(y) \; dy.
\end{split}
\end{equation}
%In the last step, we assume $U_i$ and $V_j$ are constants within $R_\omega$.
We denote
\begin{equation}
\begin{split}
 \alpha_{ij}= \int_{R_\omega} A(y)\phi_i(y) \phi_j(y) \; dy.
\end{split}
\end{equation}
It can be shown that
\begin{equation}
\begin{split}
 \alpha_{ij}  =\delta_{ij} C_i  \int_{R_\omega}\psi_j.
\end{split}
\end{equation}
From the above, we see that the macroscopic equation has the form
\[
\alpha_{ij} U_i = b_j,
\]
where 
\[
b_j=\int_{R_\omega} f \phi_j.
\]
Taking into account that $\alpha$ is a diagonal matrix,
we have
\[
U_i = {b_i\over C_i \int_{R_\omega}\psi_i}.
\]

We note that, in single continua homogenization, we obtain
\[
U_1 = b_1{1\over |R_\omega|} \int_{R_\omega} A^{-1}.
\]

\section{General case. A formal derivation.}

In this section, we present a formal derivation of generalized multicontinuum homogenization. 
The derivation makes several assumptions, which may or may not hold depending on particular problems.
We will make these assumptions as we go along.

We consider a general linear system given by
\begin{equation}
\label{eq:main0}
\begin{split}
{\partial u \over \partial t}  + Au = f, \ \text{in}\ D,
\end{split}
\end{equation}
where $A$ is a differential operator, $u$ is a vector valued solution and $D$ is the domain. 
The problem (\ref{eq:main0}) is supplemented with some appropriate initial and boundary conditions.
We next present several examples.

{\bf Example 1.} In the scalar case, 
$Au = -div(\kappa(x)\nabla u)$, where $\kappa(x)$ is a multiscale and high-contrast
coefficient.

{\bf Example 2.}  In a vector case, one can consider the elasticity problem with
$Au=-\nabla_i C_{ijkl}(x) e_{kl}(u)$, where $C_{ijkl}(x)$'s represent heterogeneous and 
high-contrast media
properties, $e_{kl}(u)=(\nabla_k u_l + \nabla_l u_k)/2$, and 
$u$ is the displacement vector. 

{\bf Example 3.} 
We can consider Example 1 in a mixed formulation as a first-order system.
In this case, $u=(p,v)$, where $p$ and $v$ solve $\kappa^{-1} v + \nabla p=0$,
$div(v)=f$.

{\bf Example 4.} One can consider the first order systems, $A u = v(x) \cdot\nabla u + a(x) u$, where
$v(x)$ and $a(x)$ are highly heterogeneous fields.

We write (\ref{eq:main0}) as a variational problem
\begin{equation}
\begin{split}
({\partial u \over \partial t},v)  + a_D(u,v) = (f,v),
\end{split}
\end{equation}
where $a_D(u,v)=\int_D (A u) v$, e.g., in Example 1, $a_D(u,v)=\int_D \kappa \nabla u\cdot  \nabla v$ (assuming zero Dirichlet boundary conditions).

In the multicontinuum homogenization, we assume that in each RVE, $R_\omega$, there exist functions
$\psi_i$ ($i$ refers to continua, $\psi_i$ can be a characteristic function of subregion), such that
\[
U_i(x_\omega^*)={\int_{R_\omega} u \psi_i \over \int_{R_\omega}\psi_i }
\]
are macroscopic variables,
where $x_\omega^*$ is a point in $R_\omega$. 
 {\bf One main assumption is that $U_i$'s are smooth functions if we consider them over all RVEs.}
Next, we present the steps in deriving macroscopic equations.

{\bf Step 1. Expansion.}

The first step consists of expanding the solution $u$ in terms of macroscopic variables.
The coefficients in front of them, denoted by $\phi_i$'s, 
represent the local microscopic solution
in RVE. 
We  consider the expansion of the solution $u$ as
\begin{equation}
\label{eq:expansion1}
\begin{split}
u = \phi_m U_m +\phi_m^j \nabla_j U_m + \phi_m^{ij} \nabla_{ij}^2 U_m + ...,
\end{split}
\end{equation}
where $\nabla_j$ refers to ${\partial \over \partial x_j}$.
In this expansion, we will discuss the functions $\phi$,
 which are defined as
the solutions of local problems in RVE, $R_\omega$.

{\bf Step 2. Cell problems.}

Next, we introduce equations for $\phi_i$'s.
These equations are written in each RVE subject to some constraints. These constraints 
are related to definitions of macroscopic variables. We use Taylor's expansion concepts
in defining the local functions such that they solve local problems with constraints that
their averages with respect to $\psi_i$ behave as constants, linear functions, 
and quadratic functions.

Our first cell problem imposes constraints to represent the constants in the average
behavior of each continua (continua $m$ in (\ref{eq:cell1}))
We consider the cell problem in oversampled regions $R_\omega^+$ that contain several $R_\omega$, denoted by $R_\omega^p$.
\begin{equation}
\label{eq:cell1}
\begin{split}
A\phi_m^i &= \Gamma_{mn}^{ijp}\psi_n^pe^j \; \text{in} \; R_\omega^+, \\
\int_{R_\omega^p} \phi_m^i \psi_n^p &= \delta_{mn} e^i \int_{R_\omega^p} \psi_n^p, \ \forall p,
\end{split}
\end{equation}
where $e^i$ is the unit vector (solution $u$ is vector valued) and
 $\psi_n^p$ is the characteristic function in $R_\omega^p$.
This cell problem corresponds to appropriate energy minimizing solution
subject to the constraints.
Here and later, by $\Gamma$, we denote the Lagrange multipliers
due to constraints. 
We denote $\phi_m$ the matrix spanned
by $\phi_m^i$ (as columns).

Our second cell problem imposes constraints to represent the 
linear functions in the average
behavior of each continua. 
\begin{equation}
\label{eq:cell2}
\begin{split}
A\phi_m^{il} &= \Gamma_{mn}^{ijpl}\psi_n^pe^j \; \text{in} \; R_\omega^+, \\
\int_{R_\omega^p}  \phi_m^{il} \psi_n^p &= \delta_{mn} e^i\int_{R_\omega^p} (x_l -c_l) \psi_n^p,\ \forall p,
\end{split}
\end{equation}
where $c_l$ (later on also) is chosen such that $\int_{R_\omega} (x_l -c_l)=0$,
where $R_\omega$ is the RVE defined in the middle of $R_\omega^+$.
Similarly,
 we denote $\phi_m^l$ the matrix spanned
by $\phi_m^{il}$ (as columns).

 We can also define higher-order cell problems. 
The next cell problem imposes constraints to represent the quadratics in the average
behavior of each continua. 
\begin{equation}
\begin{split}
A\phi_m^{ijl} &= \Gamma_{mn}^{ijlps}\psi_n^pe^s \; \text{in} \; R_\omega^+, \\
\int_{R_\omega^p} \phi_m^{ijl} \psi_n^p &= \delta_{mn} e^i\int_{R_\omega^p} (x_jx_l -c_{jl}) \psi_n^p,\ \forall p.
\end{split}
\end{equation}
Similarly,
 we denote $\phi_m^{lj}$ the matrix spanned
by $\phi_m^{ijl}$ (as columns).

We note that in our cell problems, we solve for each component of
the vector solutions. In some applications, one
can lump some components
if some relations between components of the vector are known apriori.

{\bf The decay of cell solutions.}
Existence and uniqueness can be shown in most cases for 
positive symmetric operators 
with appropriate norms. In general, we need inf-sup condition for 
well-posedness of cell
problems \cite{brezzi1974existence}.
We note that the decay of local solution away from $R_\omega$ (middle RVE)
 is important. The latter indicates
a correct computation of the local problems. In some cases, 
one can use appropriate
local boundary conditions when the information is available
about the global solution.

{\bf Step 3. Substitution in the variational formulation.}

In this step, we use $u$ and $v$ expansion in the fine-grid formulation of the problem.
In particular, we have 
\begin{equation}
\begin{split}
u = \phi_m^i U_m^i + \phi_m^{il} \nabla_l U_m^i + \phi_m^{ilp} \nabla^2_{lp} U_m^i\\
v = \phi_m^i V_m^i + \phi_m^{il} \nabla_l V_m^i+ \phi_m^{ilp} \nabla^2_{lp} V_m^i. 
\end{split}
\end{equation}

We substitute and get the following equation (we use matrix notations
for $\phi$'s)
\begin{equation}
\label{eq:expan1}
\begin{split}
& \int_D {\partial \over \partial t}(\phi_m U_m + \phi_m^{l} \nabla_l U_m + \phi_m^{lp} \nabla^2_{lp} U_m)
(\phi_n V_n + \phi_n^k \nabla_k V_n+\phi_n^{ks} \nabla_{ks}^2 V_n)\\
& a_D(\phi_m U_m, \phi_n V_n) + 
a_D (\phi_m U_m, \phi_n^k \nabla_k V_n) + 
a_D (\phi_m U_m, \phi_n^{ks} \nabla_{ks}^2 V_n) + \\
& a_D(\phi_m^l \nabla_l U_m, \phi_n V_n) + 
a_D (\phi_m^l \nabla_l U_m, \phi_n^k \nabla_k V_n) + 
a_D (\phi_m^l\nabla_l U_m, \phi_n^{ks} \nabla_{ks}^2 V_n) + \\
& a_D(\phi_m^{lp} \nabla_{lp}^2 U_m, \phi_n V_n) + 
a_D (\phi_m^{lp} \nabla_{lp}^2 U_m, \phi_n^k \nabla_k V_n) + 
a_D(\phi_m^{lp}\nabla_{lp}^2 U_m,\phi_n^{ks} \nabla_{ks}^2 V_n) \\
=
&\int_D f(\phi_n V_n + \phi_n^k \nabla_k V_n+\phi_n^{ks} \nabla_{ks}^2 V_n).
\end{split}
\end{equation}
Our next two steps include using RVE concepts and taking into account that $U_i$ and $V_i$ are smooth
functions.

{\bf Step 4. Integral localization.}

Our next step includes dividing the integral over the coarse partition
$\omega$ and then using the RVE concept. More precisely, for each
integral and a smooth function $F$, we have
\begin{equation}
\begin{split}
\int_D F  =\sum_\omega F  \approx \sum_\omega {|\omega|\over |R_\omega|} \int_{R_\omega}F.
\end{split}
\end{equation}

{\bf Step 5.  Piecesmooth approximation of macroscopic terms.}

In this step, the macroscopic terms, $U_i$ and $V_i$ assumed to be
smooth functions and the operator $A$  acts only on 
cell problem solutions.  As before (in zero-order equation case),
we take the macroscopic variables out of the integrals over $R_\omega$.
To demonstrate this step, we consider only two term expansion in (\ref{eq:expansion1}) writing the integrals over RVE.

More precisely, the terms in the equation
(\ref{eq:expan1}) have the following forms in $R_\omega$.
\begin{equation}
\label{eq:expan2}
\begin{split}
&\sum_\omega{|\omega|\over |R_\omega|} (\int_{R_\omega} \phi_m\phi_n){\partial \over \partial t} U_m V_n + \sum_\omega{|\omega|\over |R_\omega|} (\int_{R_\omega} \phi_m\phi_n^k){\partial \over \partial t} U_m \nabla_k V_n + \\
&\sum_\omega{|\omega|\over |R_\omega|} (\int_{R_\omega} \phi_m^l\phi_n){\partial \over \partial t} \nabla_l U_m V_n + \sum_\omega{|\omega|\over |R_\omega|} (\int_{R_\omega} \phi_m^{l}\phi_n^k){\partial \over \partial t} \nabla_l U_m \nabla_k V_n + \\
&+
\sum_\omega{|\omega|\over |R_\omega|} a_{R_\omega}(\phi_m,\phi_n) U_mV_n + 
\sum_\omega{|\omega|\over |R_\omega|} a_{R_\omega} (\phi_m, \phi_n^k ) U_m\nabla_k V_n+ \\
&\sum_\omega{|\omega|\over |R_\omega|} a_{R_\omega}(\phi_m^l , \phi_n ) \nabla_l U_m V_n+ 
\sum_\omega{|\omega|\over |R_\omega|} a_{R_\omega} (\phi_m^l, \phi_n^k ) \nabla_l U_m\nabla_k V_n + \\
=&\sum_\omega{|\omega|\over |R_\omega|} (\int_{R_\omega} f\phi_n) V_n + \sum_\omega{|\omega|\over |R_\omega|} (\int_{R_\omega}f\phi_n^k)\nabla_k V_n.
%\int_D f(\phi_n V_n + \phi_n^k \nabla_k V_n+\phi_n^{ks} \nabla_{ks}^2 V_n)
\end{split}
\end{equation}

In Equation (\ref{eq:expan2}), 
we further take into account that $U_i$ and $V_i$ are smooth
functions defined in $D$ and get the following macroscopic 
equation for $U_i$ (in strong form)
\begin{equation}
\begin{split}
A_{nm}{\partial \over \partial t} U_m+
B_{nm}U_m + B_{nm}^i\nabla_iU_m - 
\nabla_k\overline{B}^k_{nm}U_m - \nabla_k(B_{nm}^{ik}\nabla_iU_m)=b_n. 
\end{split}
\end{equation}
Here, we neglect the second, third, and fourth terms in (\ref{eq:expan2}).
The latter is because $\phi_n^k$ is of order RVE size, while
$\phi_n$ is of order $O(1)$, in general. Because the coefficients
in the operator $A$ have high-contrast properties, 
we can not neglect these terms. We will remark on this later.
The coefficients $A$'s and $B$'s are defined from (\ref{eq:expan2}). More precisely,
\begin{equation}
\begin{split}
&A_{nm}={1 \over |R_\omega|}\int_{R_\omega} \phi_m\phi_n, \ \ b_n ={1 \over |R_\omega|}\int_{R_\omega}f\phi_n^k,\\
&B_{nm}= {1 \over |R_\omega|}a_{R_\omega}(\phi_m,\phi_n),\ \ B_{nm}^i= {1 \over |R_\omega|}a_{R_\omega}(\phi^i_m,\phi_n),\\
&\overline{B}^k_{nm}={1 \over |R_\omega|}a_{R_\omega}(\phi_m,\phi_n^k), \ \ 
B_{nm}^{ik}= {1 \over |R_\omega|}a_{R_\omega}(\phi^i_m,\phi^k_n).\\
\end{split}
\end{equation}

If we use the second-order expansion, the macroscopic equation will have the following form
\begin{equation}
\begin{split}
&A_{nm}{\partial \over \partial t} U_m+\\
&B_{nm}U_m + B_{nm}^i\nabla_iU_m + B_{nm}^{ij}\nabla^2_{ij}U_m-\\
&\nabla_k(B^k_{nm}U_m) - \nabla_k(B_{nm}^{ik}\nabla_iU_m) - \nabla_k (B_{nm}^{ijk}\nabla^2_{ij}U_m)+\\
&\nabla^2_{kp}(B^{kp}_{nm}U_m) + \nabla^2_{kp}(B_{nm}^{ikp}\nabla_iU_m) + \nabla^2_{kp}(B_{nm}^{ijkp}\nabla^2_{ij}U_m)=b_n.
\end{split}
\end{equation}

Next, we make several remarks.

First, different terms in the macroscopic equation can have negligible weights.
In general, $\phi_k$'s (the cell solutions accounting for the averages)
are of order $O(1)$, while
 $\phi_k^n$'s 
(the cell solutions accounting for the gradients)
are of order $O(\epsilon)$, where $\epsilon$ is the RVE size
(see \cite{efendiev2023multicontinuum}). For this reason, we
have neglected some terms in the time derivative terms and source terms.
However, because of high-contrast coefficients, one can not neglect
different terms that stem from  $\phi_k$'s or from  $\phi_k^n$'s.
In \cite{efendiev2023multicontinuum}, 
we show that the zero-order terms are important 
when there is high contrast. More precisely, the reaction terms
scale as the inverse of the RVE size. If the effective diffusivity is high,
then the reaction and diffusion terms balance each other. Otherwise, 
one can show that there is no multicontinuum and our macroscopic equations
result to single continuum homogenization.

Our second remark is regarding the definition of the continua.
Throughout the paper, we assume that $\psi_i$'s are associated with 
subregions defined apriori. In general, one can use  
spatial functions for $\psi_i$, for example, defined via local
spectral problems as it is done in nonlocal multicontinua approach
or GMsFEM \cite{GMsFEM13,vasilyeva2019nonlocal}.

\subsection{Example. A scalar elliptic equation}

This example is discussed in \cite{efendiev2023multicontinuum}. We briefly
mention it here. We will focus on multicontinuum expansion,
macroscopic equations, and constraints, for simplicity, and do not write down
the cell problem equations (cf. \ref{eq:cell1}).
The multicontinuum expansion is $u=\phi_i U_i+\phi^m_{i}\nabla_m U_i$, where
cell solutions have constraints for $\phi_m$
\begin{equation}
\label{eq:cell_ex1}
\begin{split}
\int_{R_\omega^p} \phi_m \psi_n^p = \delta_{mn} \int_{R_\omega^p} \psi_n^p
\end{split}
\end{equation}
and for $\phi_m^l$
\begin{equation}
\label{eq:cell_ex12}
\begin{split}
\int_{R_\omega^p}  \phi_m^{l} \psi_n^p  = \delta_{mn}\int_{R_\omega^p} (x_l -c_l) \psi_n^p.
\end{split}
\end{equation}
Note that the equations for $\phi_m$ and for $\phi_m^{l}$ are solved
separately.

The macroscopic equations have the following form
\begin{equation}
\label{eq:macro11}
\begin{split}
B_{nm}U_m + B_{nm}^i\nabla_iU_m -
\nabla_k(\overline{B}^k_{nm}U_m) - \nabla_k(B_{nm}^{ik}\nabla_iU_m) =b_n.
\end{split}
\end{equation}

\subsection{Example. A system of elliptic equations}
We consider 
\begin{equation}
-{\partial \over \partial x_k} (A^{kl}_{ji}
{\partial \over \partial x_l} u_i) = f_j.
\end{equation}
The multicontinuum expansion has the following form
\[
u_i = \phi_{mij}U_{mj} + \phi_{mij}^k\nabla_k U_{mj},
\]
where the cell problems have the constraints
for $\phi_{imn}$
\begin{equation}
\label{eq:cell_ex21}
\begin{split}
\int_{R_\omega^p} \phi_{mij} \psi_n^p = \delta_{mn}\delta_{ij} \int_{R_\omega^p} \psi_n^p
\end{split}
\end{equation}
and for $\phi_{mij}^{l}$
\begin{equation}
\label{eq:cell_ex22}
\begin{split}
\int_{R_\omega^p}  \phi_{mij}^{l} \psi_n^p = \delta_{mn}\delta_{ij}\int_{R_\omega^p} (x_l -c_l) \psi_n^p .
\end{split}
\end{equation}

For example, for two equations, we have
\begin{equation}
\begin{bmatrix}
    u_1 \\
    u_2\\
  \end{bmatrix}
=
\begin{bmatrix}
    \phi_{j11} &  \phi_{j12}   \\
    \phi_{j21} & \phi_{j22} \\
  \end{bmatrix}
\begin{bmatrix}
    U_{j1} \\
     U_{j2}\\
\end{bmatrix}
+
\begin{bmatrix}
   \phi^m_{j11} &\phi^m_{j12} \\
   \phi^m_{j21} &  \phi^m_{j22} \\
  \end{bmatrix}
\begin{bmatrix}
    \nabla_m U_{j1}\\
\nabla_m U_{j2}
  \end{bmatrix}.
\end{equation}
The constraints are the following
\begin{equation}
\label{eq:cell_ex23}
\begin{split}
\int_{R_\omega^p} 
\begin{bmatrix}
    \phi_{m11}    \\
    \phi_{m21} \\
  \end{bmatrix}
 \psi_n^p = \delta_{mn}
\begin{bmatrix}
    1  \\
    0\\
  \end{bmatrix}
 \int_{R_\omega^p} \psi_n^p,
\ \ \ 
\int_{R_\omega^p} 
\begin{bmatrix}
    \phi_{m12}    \\
    \phi_{m22} \\
  \end{bmatrix}
 \psi_n^p = \delta_{mn}
\begin{bmatrix}
    0  \\
    1\\
  \end{bmatrix}
 \int_{R_\omega^p} \psi_n^p,
\end{split}
\end{equation}
\begin{equation}
\label{eq:cell_ex231}
\begin{split}
\int_{R_\omega^p} 
\begin{bmatrix}
    \phi^k_{m11}    \\
    \phi^k_{m21} \\
  \end{bmatrix}
 \psi_n^p = \delta_{mn}
\begin{bmatrix}
    1  \\
    0\\
  \end{bmatrix}
 \int_{R_\omega^p} (x_k-c_k)\psi_n^p,
\\\ 
\int_{R_\omega^p} 
\begin{bmatrix}
    \phi^k_{m12}    \\
    \phi^k_{m22} \\
  \end{bmatrix}
 \psi_n^p = \delta_{mn}
\begin{bmatrix}
    0  \\
    1\\
  \end{bmatrix}
 \int_{R_\omega^p} (x_k-c_k)\psi_n^p.
\end{split}
\end{equation}

The macroscopic equations have the following form
\begin{equation}
\begin{split}
B_{ijnm}U_{jm} + B_{ijnm}^l\nabla_lU_{jm} -
\nabla_k(B^k_{ijnm}U_{jm}) - \nabla_k(B_{ijnm}^{lk}\nabla_lU_{jm}) =b_{in}.
\end{split}
\end{equation}
It can be shown that the second and third terms cancel each other and
the scaling of 
$B_{nm}$ is of order $1/\epsilon^2$, where $\epsilon$ is RVE size. 
Because of high contrast, this term can balance with the diffusion term.

\subsection{Mixture theory and its relation}

Here, we briefly note that one can also 
derive general multicontinuum equations using mixture theory
\cite{rajagopal1995mechanics,truesdell1984thermodynamics,malek20};
however, precise micro and macro relations can not be derived from this theory.
Mixture theory specifies several model classes \cite{malek20}. 
One that is suitable for our models is Class II, where
 $N$ balances of mass for N components of the mixture and also $N$ balances of linear momentum for N components of mixture are formulated. 
In this case, the equations have the following form
\begin{equation}
\begin{split}
{\partial \rho_i\over \partial t} + div (\rho_i  v_i) = m_i,\ \ \sum_i m_i =0,\\
{\partial \rho_i v_i \over \partial t} + div (\rho_i v_i \otimes v_i)=div(\mathcal{T}_i) + \mathcal{Q}_i+ m_i v_i,\ \ \sum_i (\mathcal{Q}_i + m_i v_i) =0.
\end{split}
\end{equation}
Here, we use a simplified formulation from 
\cite{malek20},  and use
the notations from \cite{malek20}, 
where $\rho_i$ is the density of $i$th component,
$v_i$ is the velocity, $m_i$ is the exchange terms for mass conservation, 
$\mathcal{T}_i$
is the stress tensor,  and $\mathcal{Q}_i$ is the exchange terms for
momentum. 

To derive a multicontinuum equations, we consider solid and two
fluid continua mixture. 
For momentum equations, we have (ignoring gravity)
\begin{equation}
\begin{split}
{\partial \rho_1^f v_1^f \over \partial t} + div (\rho_1^f v_1^f \otimes v_1^f)=div(\mathcal{T}^f_1) + \mathcal{Q}^f_1,\\
{\partial \rho_2^f v_2^f \over \partial t} + div (\rho_2^f v_2^f \otimes v_2^f)=div(\mathcal{T}^f_2) + \mathcal{Q}^f_2,\\\
{\partial \rho^s v^s \over \partial t} + div (\rho^s v^s \otimes v^s)=div(\mathcal{T}_s) +\mathcal{Q}^s,
\end{split}
\end{equation}
where $\mathcal{Q}^s=-\mathcal{Q}^f_1-\mathcal{Q}^f_2$,
$s$ denotes the solid and $f$ denotes the fluid.
 It is assumed that
$v^s \approx 0$, $\mathcal{T}^f_i=p_i I$, $i=1,2$,  
$\mathcal{Q}^f_i=\kappa_i^{-1} v_i$,
and
the flow is steady-state and slow. In the mass conservation equations,
\begin{equation}
\begin{split}
{\partial \rho_1^f\over \partial t} + div (\rho_1^f v_1^f )= m^f_1,\\
{\partial \rho_2^f \over \partial t} + div (\rho_2^f v_2^f)=m^f_2,\\\
{\partial \rho^s \over \partial t} + div (\rho^s v^s)=m^s.
\end{split}
\end{equation}
We have $v^s\approx 0$, $m^s \approx 0$, and take
\begin{equation}
\begin{split}
m^f_1 = \alpha \rho_1^f (p_2-p_1)\\
m^f_2 = \alpha \rho_2^f (p_1-p_2).\\
\end{split}
\end{equation}
The resulting equations have the form of multicontinuum equations
(\ref{eq:macro11}).

\section{First-order mixed system}

We consider a first-order mixed system as an example of a system, 
where the variables are
coupled.
\begin{equation}
\label{eq:mixed_main}
\begin{split}
\kappa^{-1} v + \nabla u = 0\\
div( v) = f.
\end{split}
\end{equation}
This equation is a non-symmetric system with the solution vector $(v,u)$ and
the operator 
\begin{equation}
A=
\begin{bmatrix}
    \kappa^{-1} & \nabla \\
    div  & 0 \\
  \end{bmatrix}.
\end{equation}
The local cell problems and constraints require special attention to achieve 
a decay property. We omit this part to numerical results. We consider 
the derivation of macroscopic equations. In general, as before, one can
use various constraints and derive various macroscopic equations.

We consider piecewise constant velocity and
piecewise linear type pressure approximations at the RVE level.
We use different notations because differing notations for variables.
In this case, we have the following expansion
\begin{equation}
\label{eq:expan11}
\begin{split}
&v_s=\phi_{is}^{vu} U_i + \phi_{ism}^{vu} \nabla_m U_i + \phi_{isk}^{vv} V_{ik}\\  
&u = \phi^{uu}_i U_i + \phi^{uu}_{im} \nabla_m U_i + \phi_{ik}^{uv} V_{ik}.\\
\end{split}
\end{equation}
Here, $i$ refers to the continua, $(\phi^{uv},\phi^{uu})$ represents 
the cell solutions
with zero constraints on $v$ and $(\phi^{vu},\phi^{vv})$ represents cell 
solutions with zero constraints on $u$ 
(see Section 5, (\ref{eq:mixed_c1})-(\ref{eq:mixed_c3})).

We multiple the mixed system (\ref{eq:mixed_main}) 
by 
\begin{equation}
\begin{bmatrix}
    \phi^{vu}_j Q_j + \phi_{jl}^{vu}\nabla_l Q_j\\
    \phi^{uu}_j Q_j + \phi_{jl}^{uu}\nabla_l Q_j
  \end{bmatrix}
\end{equation}
and sum up the equations (use vector notations for simplicity)
\begin{equation}
\begin{split}
&\int_{R_\omega} (\phi^{vu}_j Q_j + \phi_{jl}^{vu}\nabla_l Q_j)\kappa^{-1}(\phi^{vu}_i U_i + \phi_{im}^{vu}\nabla_m U_i+\phi_{i}^{vv} V_{i}) \\
&+\int_{R_\omega} (\phi^{vu}_j Q_j + \phi_{jl}^{vu}\nabla_l Q_j) \nabla(\phi^{uu}_i U_i + \phi_{im}^{uu}\nabla_m U_i+\phi_{i}^{uv} V_{i})\\
&+
\int_{R_\omega} div(\phi^{vu}_i U_i + \phi_{im}^{vu}\nabla_m U_i+\phi_{i}^{vv} V_{i})(\phi^{uu}_j Q_j + \phi_{jl}^{uu}\nabla_l Q_j) \\
=&\int_{R_\omega} f (\phi^{uu}_j Q_j + \phi_{jl}^{uu}\nabla_l Q_j).
\end{split}
\end{equation}
In the global form, the equation has the form
\begin{equation}
\label{eq:mix21}
\begin{split}
\alpha^{u}_{ij}U_i + \alpha^{u}_{ijm}\nabla_m U_i  -  \nabla_m(\overline{\alpha}^{u}_{ijm} U_i)-\nabla_n(\alpha^{u}_{ijnm}\nabla_m U_i)+\beta^{u}_{ji}  V_{i}  + \beta^{u}_{jim} \nabla_m V_{i} = f^u_j.
\end{split}
\end{equation}
In our numerical simulations, we observe that the sum of two convection terms
(the second and third terms in the equation)
is small and can be neglected.

Next, we multiply the system (\ref{eq:mixed_main}) by 
\begin{equation}
\begin{bmatrix}
    \phi_{j}^{vv} Q_{j}\\
    \phi_{j}^{uv} Q_{j} 
  \end{bmatrix}
\end{equation} 
and sum up (use vector notations for simplicity)
\begin{equation}
\begin{split}
&\int_{R_\omega} \phi_{j}^{vv} Q_{j}\kappa^{-1}(\phi^{vu}_i U_i + \phi_{im}^{vu}\nabla_m U_i+\phi_{i}^{vv} V_{i}) \\
&+\int_{R_\omega} \phi_{j}^{vv} Q_{j}\nabla(\phi^{uu}_i U_i + \phi_{im}^{uu}\nabla_m U_i+\phi_{i}^{uv} V_{i})\\
&+
\int_{R_\omega} div(\phi^{vu}_i U_i + \phi_{im}^{vu}\nabla_m U_i+\phi_{i}^{vv} V_{i})\phi_{j}^{uv} Q_{j} \\
=&\int_{R_\omega} f \phi_{j}^{uv} Q_{j}.
\end{split}
\end{equation}
In the global form, the equation has the form
\begin{equation}
\label{eq:mix22}
\begin{split}
\alpha^{v}_{ij}U_i + \alpha^{v}_{ijm}\nabla_m U_i  + \beta^{v}_{ji}  V_{i} = f^v_j.
\end{split}
\end{equation}
In our numerical simulations, we observe that $\alpha^{v}_{ij}$ is small and
can be neglected.

In general, one can choose a more general representation of the velocity via
piecewise linear functions and obtain general models with higher
order derivatives.

 Note that the polynomial constraints in
the  approximation of velocity and pressures   in 
(\ref{eq:expan11}) is for homogenization and is not related to
stable polynomial approximation in finite element methods.

\section{Numerical example}

In this section, we will present some numerical examples to demonstrate 
the performance of the method in a mixed formulation. 
As we mentioned earlier that for nonsymmetric problems, it is challenging
to guarantee the decay of local solutions. Here, we propose
the following local problems for velocity $v$ and pressure $u$ in the equation
\begin{equation}
\begin{split}
\kappa^{-1} v + \nabla u &=0\\
div(v)&=q. 
\end{split}
\end{equation}

Next, we describe the local solutions for velocity and pressure. 
 We will only 
write down the constraints, the formulation of the local problem
follows from equations (\ref{eq:cell1}) and (\ref{eq:cell2}).
For the velocity constraints, we impose an intermediate domain
$R_\omega^V$, where $R_\omega^V$ is a subset of $R_\omega^+$ and contains
$R_\omega$. Moreover, we assume that $R_\omega^V$ consists of $R_\omega^p$,
where $p$ is a numeration of local domains, one of them being $R_\omega$.
We remind that the local solution has the following matrix form.
\begin{equation}
\begin{bmatrix}
    u \\
    v_s\\
  \end{bmatrix}
=
\begin{bmatrix}
    \phi^{uu}_i & \phi^{uu}_{im} & \phi^{uv}_{ik} \\ 
    \phi^{vu}_{is} & \phi^{vu}_{ism} & \phi^{vv}_{isk} \\ 
  \end{bmatrix}
\begin{bmatrix}
    U_i \\
    \nabla_m U_i\\
      V_{ik}
  \end{bmatrix}
\end{equation}
The local constraints for $\phi$'s are imposed column by column.
The constraints
are the following
\begin{equation}
\label{eq:mixed_c1}
\begin{split}
\int_{R_\omega^p}\phi^{uu}_i \psi_j &= \delta_{ij} \int_{R_\omega^p}\psi_j, \ \forall R_\omega^p\subset R_\omega^+,\\
\int_{R_\omega^p}\phi^{vu}_{is} \psi_j &= 0, \ \forall R_\omega^p\subset R_\omega^V,\\
\end{split}
\end{equation}
and
\begin{equation}
\label{eq:mixed_c2}
\begin{split}
\int_{R_\omega^p}\phi^{uu}_{im} \psi_j &= \delta_{ij} \int_{R_\omega^p}(x_m -c_m)\psi_j, \ \forall R_\omega^p\subset R_\omega^+,\\
\int_{R_\omega^p}\phi^{vu}_{ism} \psi_j &= 0, \ \forall R_\omega^p\subset R_\omega^V,\\
\end{split}
\end{equation}
and
\begin{equation}
\label{eq:mixed_c3}
\begin{split}
\int_{R_\omega^p}\phi^{uv}_{ik} \psi_j &= 0, \ \forall R_\omega^p\subset R_\omega^+,\\
\int_{R_\omega^p} \phi^{vv}_{isk} \psi_j &= \delta_{ij} \delta_{sk} \int_{R_\omega^p} \psi_j, \ \forall R_\omega^p\subset R_\omega^V.
\end{split}
\end{equation}

In the calculations of macroscopic domains, we use another intermediate
domain $R_\omega^I$, which is a subset of $R_\omega^+$ and contains
$R_\omega^V$. The local expansion is given by (\ref{eq:expan11}).

\begin{figure}[h!]
\centering
\includegraphics[scale=0.3]{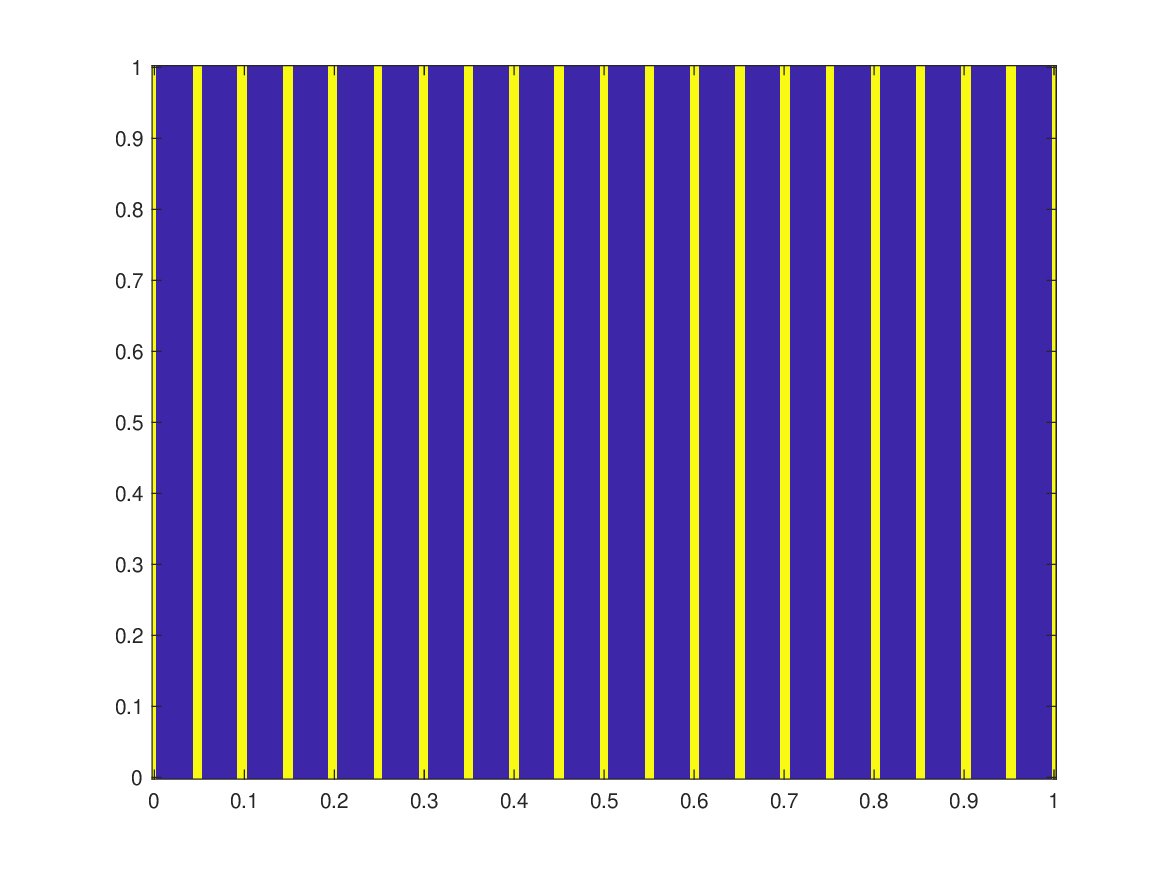} 
\includegraphics[scale=0.3]{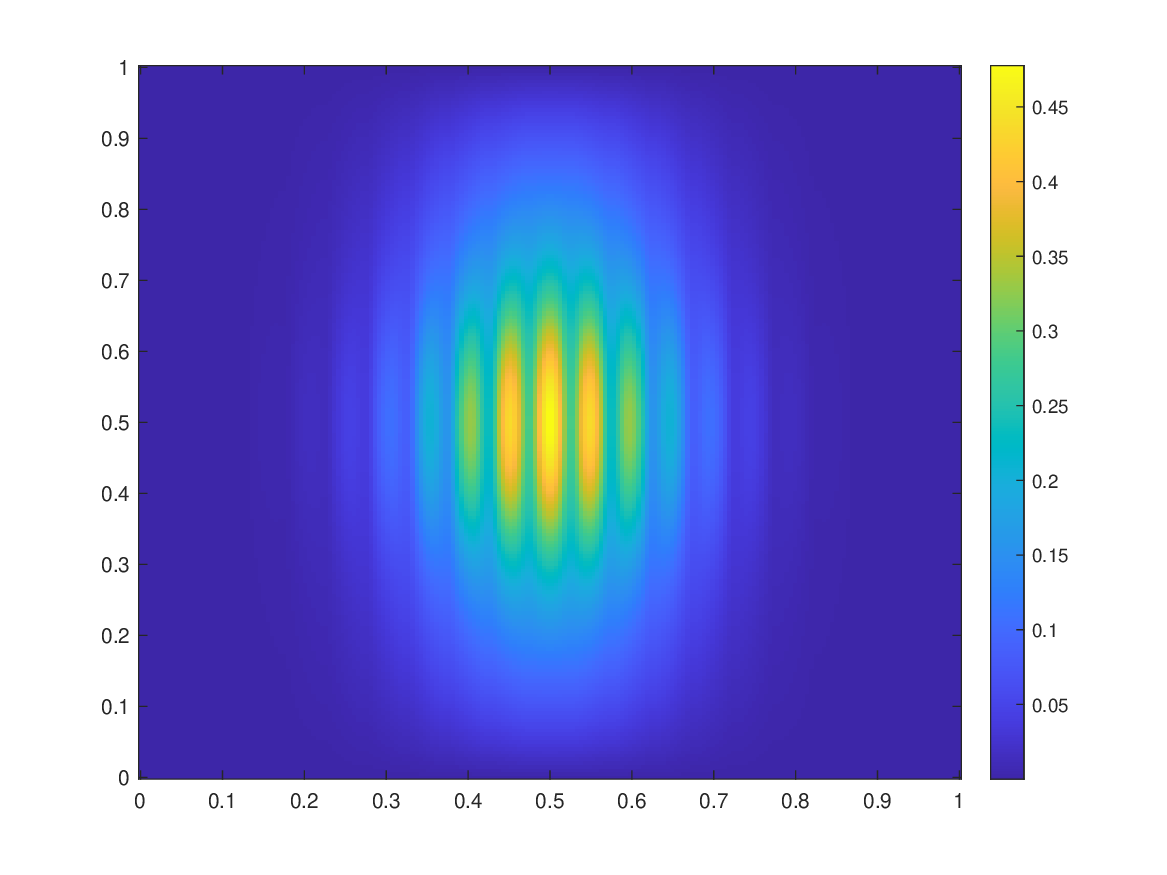}
\caption{Case 1. Left: Parameter $\kappa$. Right: Reference solution.}
\label{fig:perm1}
\end{figure}

In the first example, we consider the layered medium 
depicted on Figure \ref{fig:perm1}.
The permeability field
$\kappa$
has a period denoted by $\epsilon$. We denote the low
conductivity region and the high conductivity region of $\kappa$
by $\Omega_{1}$ and $\Omega_{2}$, respectively.
The source term
$f$ and conductivity $\kappa$ are as follows
\[
f(x)=\begin{cases}
1000\min\{\kappa\}e^{-40|(x-0.5)^{2}+(y-0.5)^{2}|} & x\in\Omega_{1}\\
e^{-40|(x-0.5)^{2}+(y-0.5)^{2}|} & x\in\Omega_{2}
\end{cases}
\]
and
\[
\kappa(x)=\begin{cases}
\cfrac{\epsilon}{10000} & x\in\Omega_{1}\\
\cfrac{1}{100\epsilon} & x\in\Omega_{2}
\end{cases}
\]

We divide the 
 computational domain $\Omega$  into $M\times M$
coarse grid. The coarse mesh size $H$ is defined as $H=1/M$. 
We consider the whole coarse grid element as an RVE for
the corresponding coarse element. The oversampling RVE $R_\omega^{+}$
(or $\omega^+$)
for each coarse RVE $\omega$ is defined as an extension of $K$ (target
coarse block) by
$l$ layers of coarse grid element, where $l$ will be changed in simulations.

We define the relative $L^{2}$- error in $\Omega_{1}$ and the
relative $L^{2}$- error in $\Omega_{2}$ by 
\[
e_{2}^{(i)}=\cfrac{\sum_{K}|\cfrac{1}{|K|}\int_{K}U_{i}-\cfrac{1}{|K\cap\Omega_{i}|}\int_{K\cap\Omega_{i}}u|^{2}}{\sum_{K}|\cfrac{1}{K\cap\Omega_{i}}\int_{K\cap\Omega_{i}}u|^{2}}.
\]
$K$ denotes the RVE, which is taken to be $\omega$.
%This represents the $L_2$ error of our proposed approach.

\begin{figure}[h!]
\centering
\includegraphics[scale=0.3]{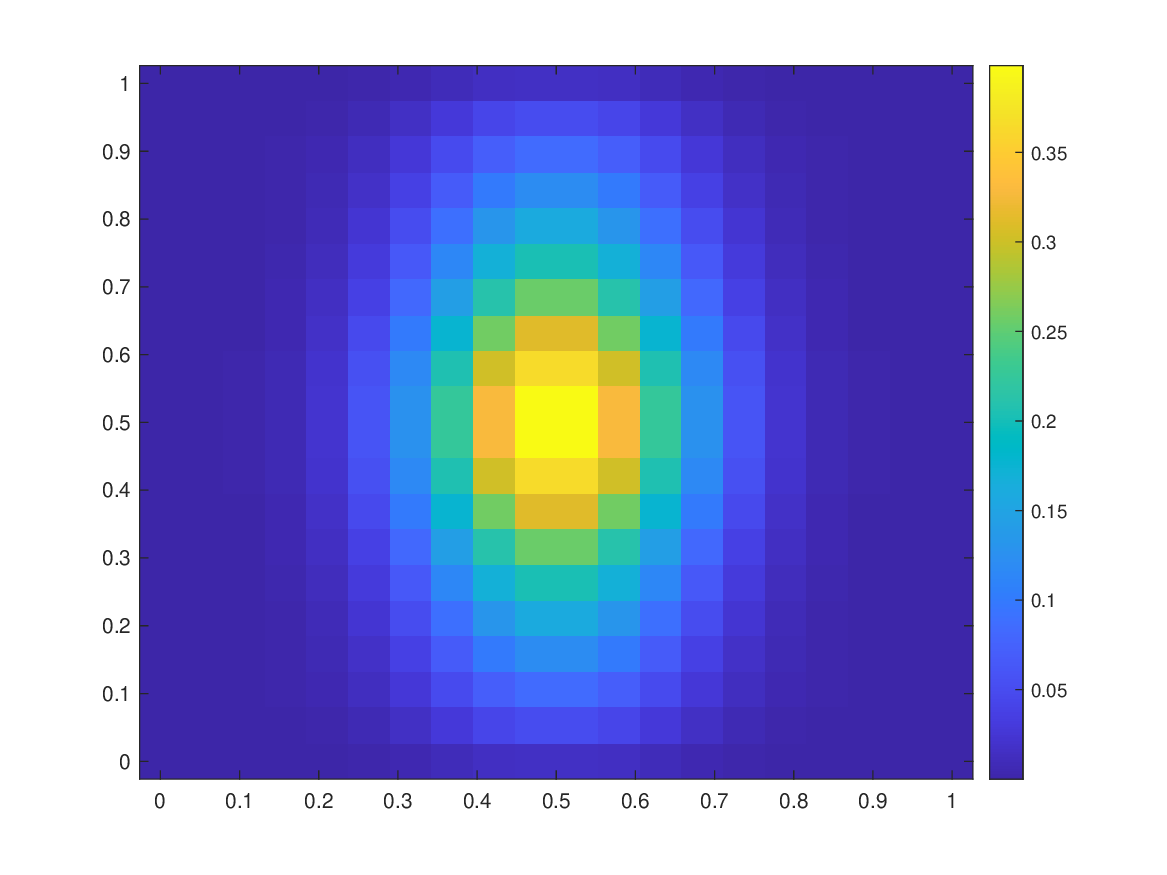} 
\includegraphics[scale=0.3]{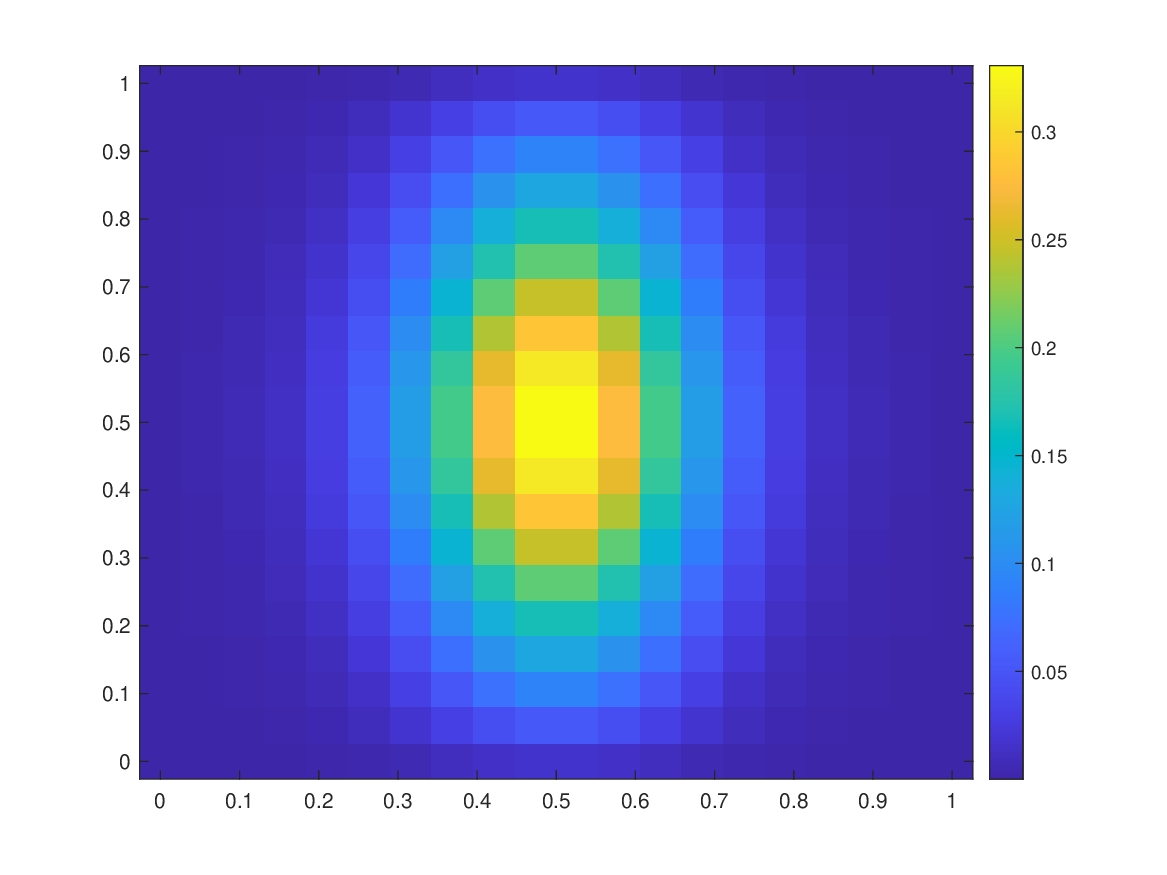}
\includegraphics[scale=0.3]{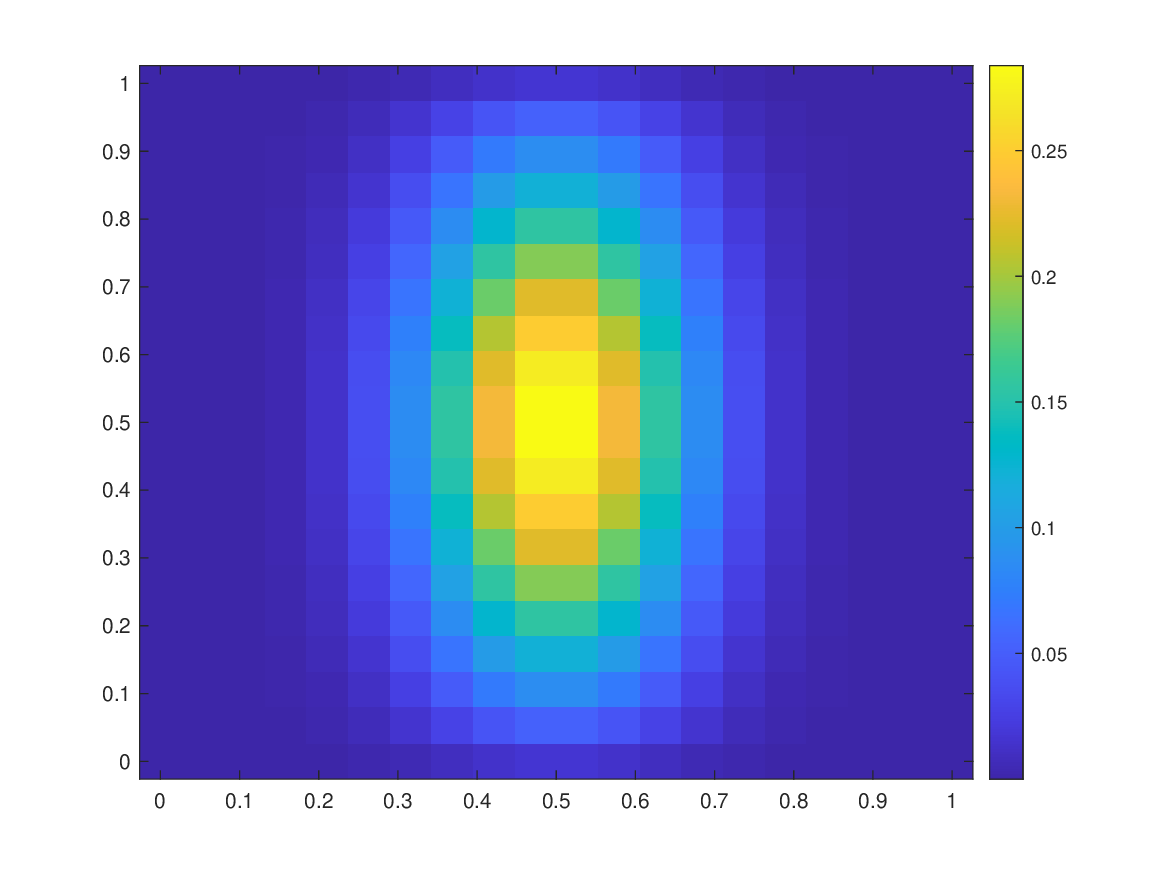} 
\includegraphics[scale=0.3]{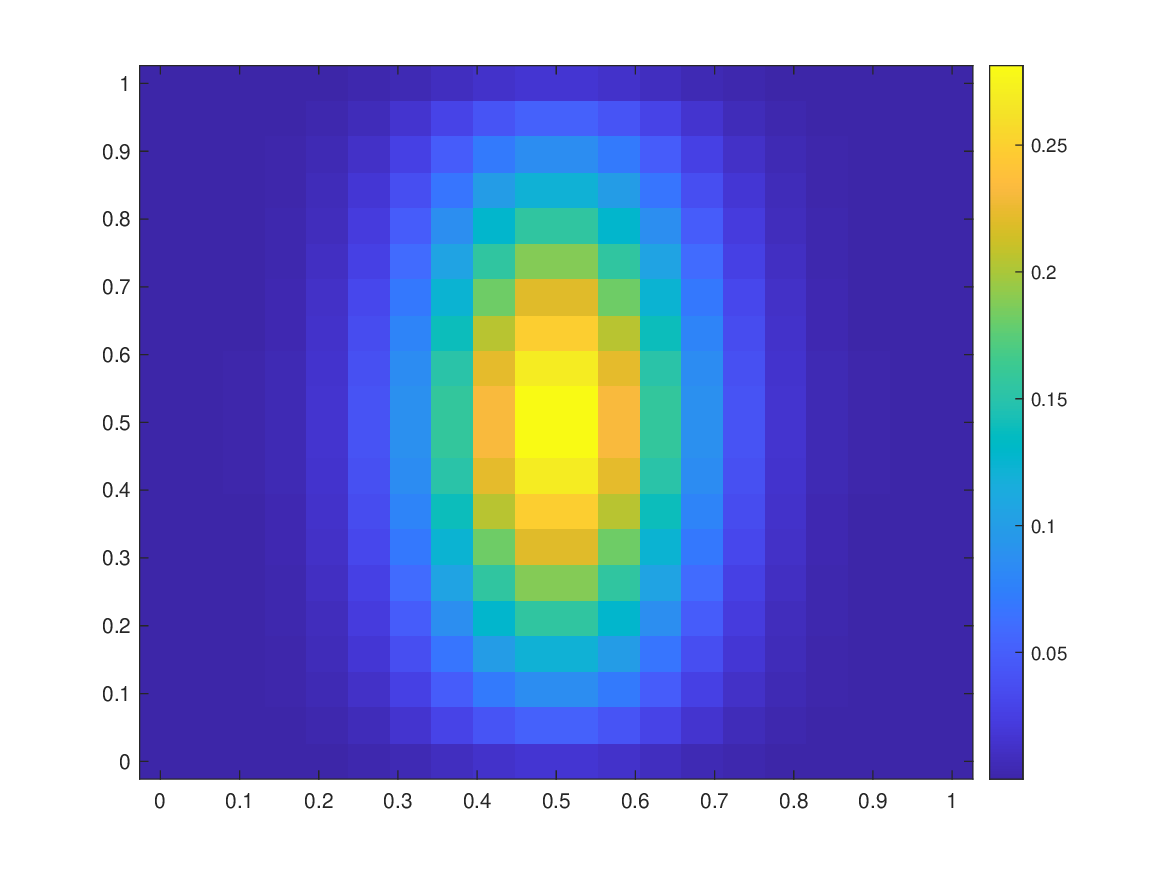}
\caption{Case 1. Top-Left: reference average solution in $\Omega_1$.  Top-Right: homogenized average solution in $\Omega_1$.  Bottom-Left: reference average solution in $\Omega_2$.  Bottom-Right: homogenized average solution in $\Omega_2$.}
\label{fig:compare_case1}
\end{figure}

For the first case, we take the fine-mesh size to be $H\epsilon$.
We present $e_2^{(i)}$ in Table \ref{tab:case1}.
 First, we observe that the proposed
approach provides an accurate approximation of the averaged solution
as we decrease the mesh size.
In Figure \ref{fig:compare_case1}, we depict upscaled solutions
and corresponding averaged fine-scale solutions. We observe that these
solutions are very close.
In the first table, we decrease the coarse-mesh size and the
period size. 
In standard numerical homogenization methods,
this gives a resonance error (stagnating errors).
Here, by choosing an appropriate number of layers, we observe that
the error remains small. In the second table,  we observe convergence as we decrease the mesh size and fix
$\epsilon$. In general,
we expect a certain threshold error 
due to fine-scale discretization, which is
used to compute the solution.

\begin{table}[h!]
\centering
\begin{tabular}{|c|c|c|c|}
\hline 
$H$ & $\epsilon$ & $e_{2}^{(1)}$ & $e_{2}^{(2)}$\tabularnewline
\hline 
$1/10$ & $1/10$ & $27.19\%$ & $6.21\%$\tabularnewline
\hline 
$1/20$ & $1/20$ & $11.63\%$ & $1.19\%$\tabularnewline
\hline 
$1/40$ & $1/40$ & $3.25\%$ & $0.88\%$\tabularnewline
\hline 
\end{tabular} 
\begin{tabular}{|c|c|c|c|}
\hline 
$H$ & $\epsilon$ & $e_{2}^{(1)}$ & $e_{2}^{(2)}$\tabularnewline
\hline 
$1/10$ & $1/10$ & $27.19\%$ & $6.21\%$\tabularnewline
\hline 
$1/10$ & $1/20$ & $12.79\%$ & $2.43\%$\tabularnewline
\hline 
$1/10$ & $1/40$ & $6.32\%$ & $1.70\%$\tabularnewline
\hline 
\end{tabular}
\caption{Error comparison for Case 1.}
\label{tab:case1}
\end{table}

\begin{figure}[h!]
\centering
\includegraphics[scale=0.3]{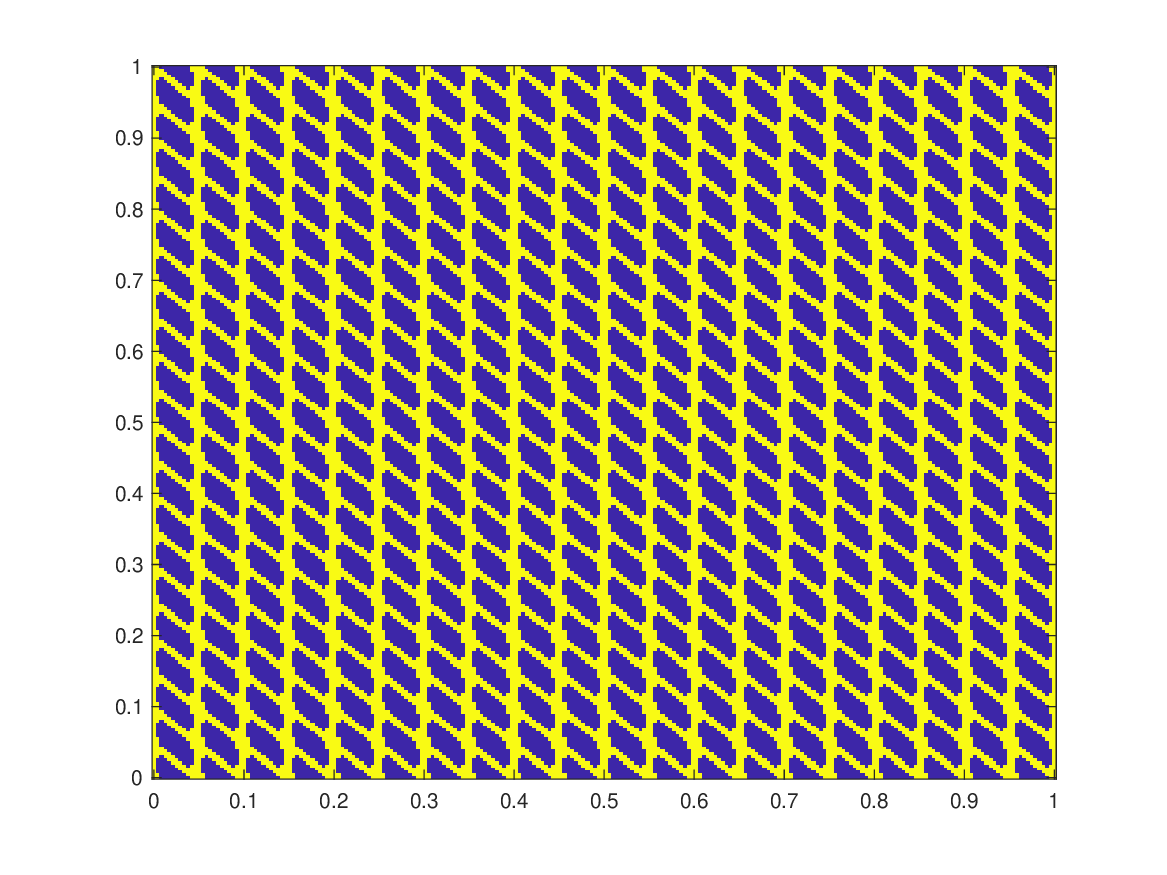} 
\includegraphics[scale=0.3]{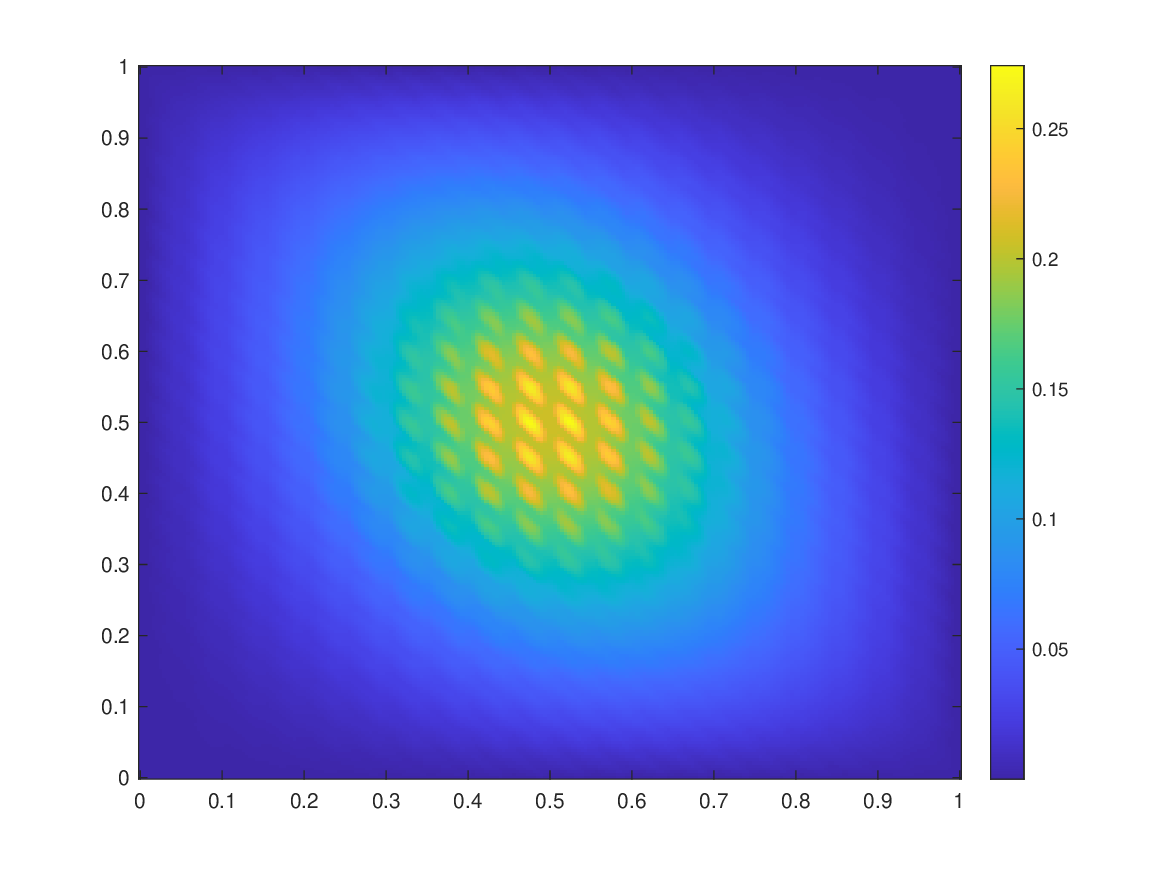}
\caption{Case 2. Left: Parameter $\kappa$. Right: reference solution.}
\label{fig:case2_perm}
\end{figure}

\begin{figure}[h!]
\centering

\includegraphics[scale=0.3]{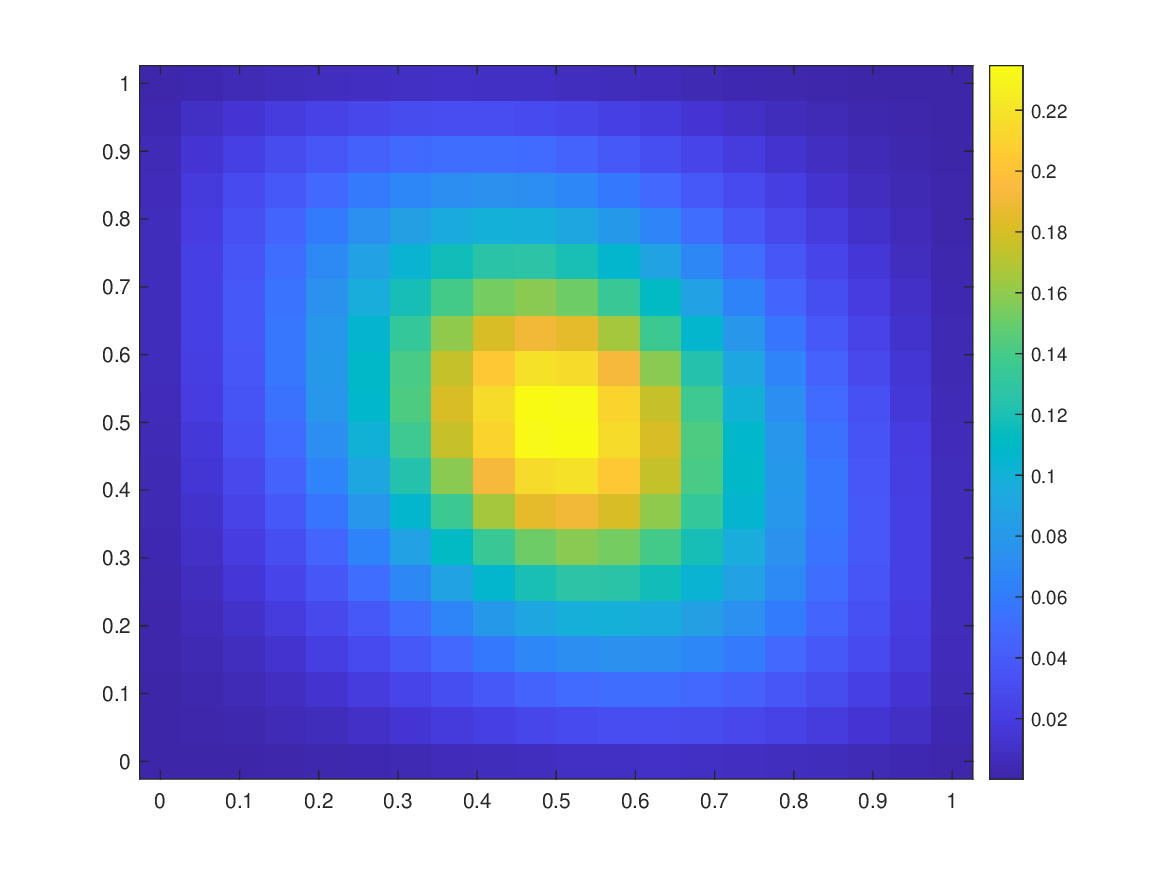} 
\includegraphics[scale=0.3]{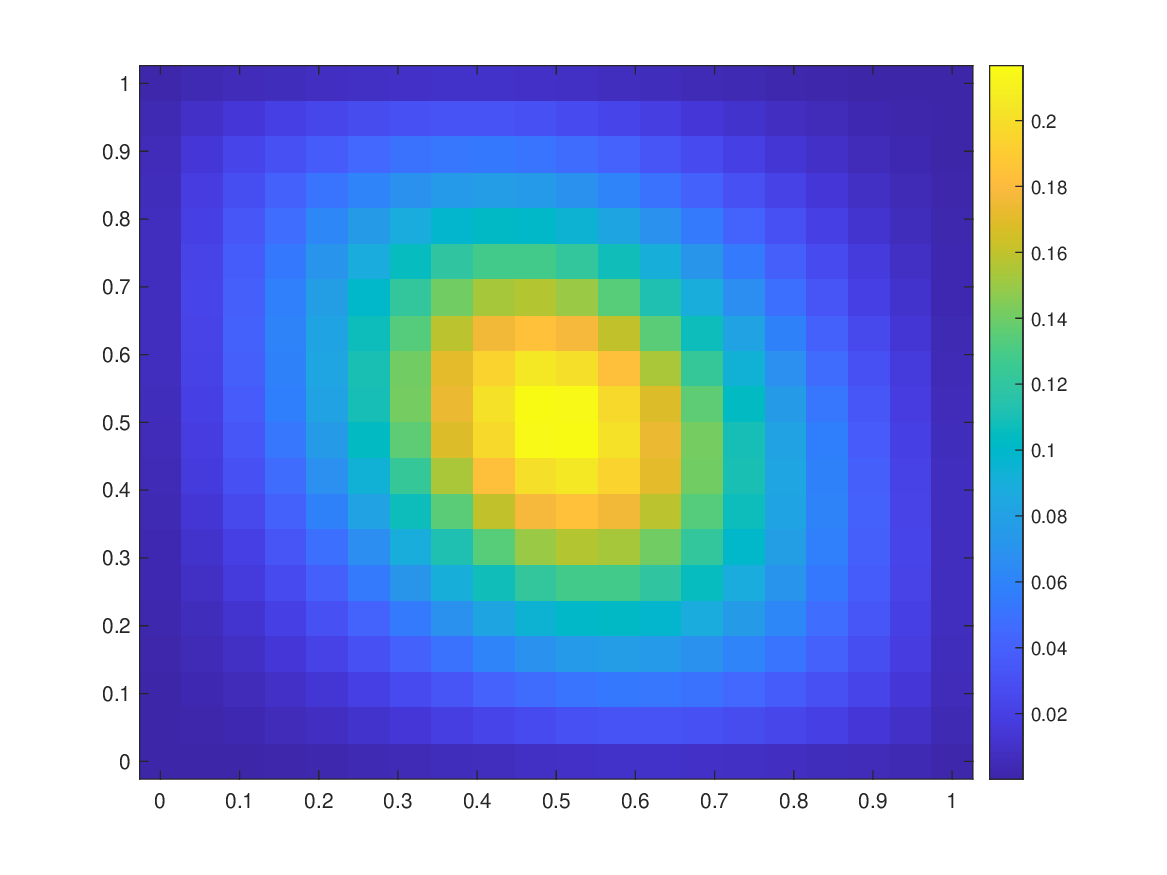}

\includegraphics[scale=0.3]{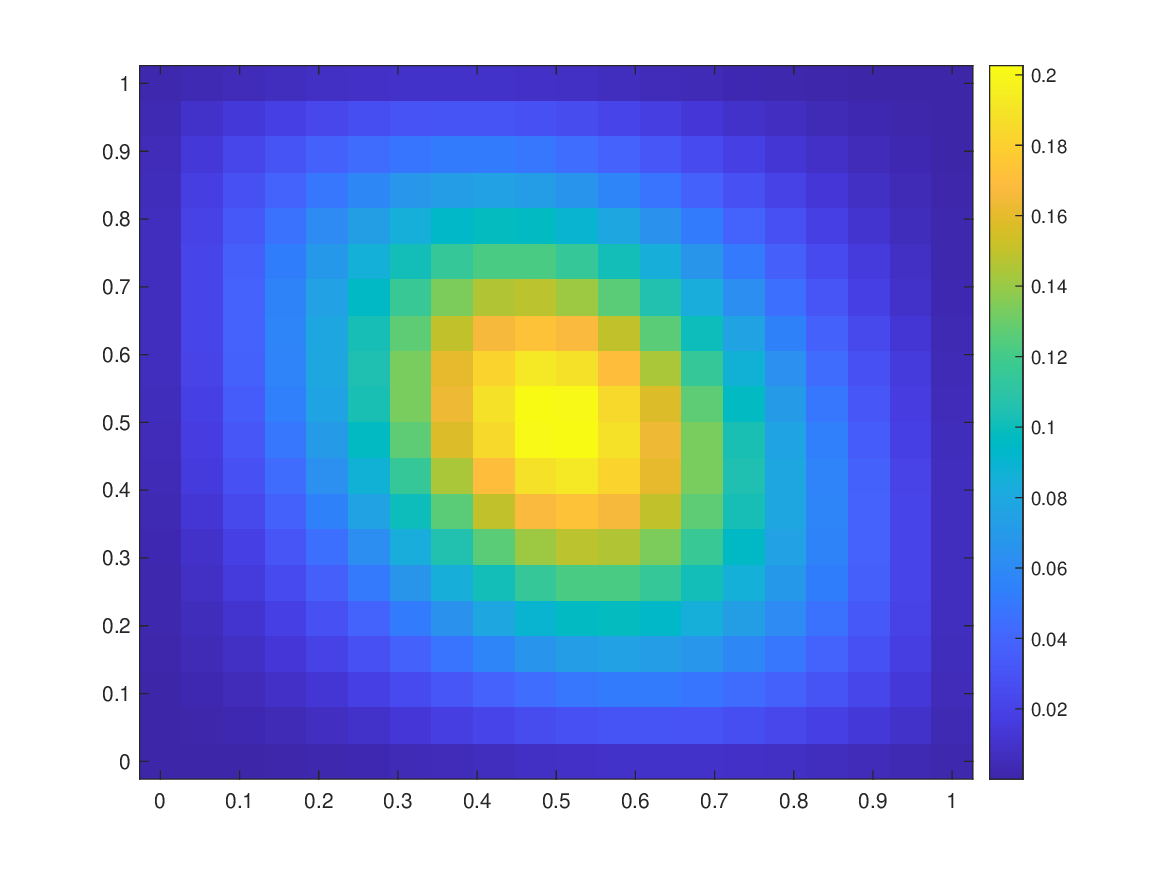} 
\includegraphics[scale=0.3]{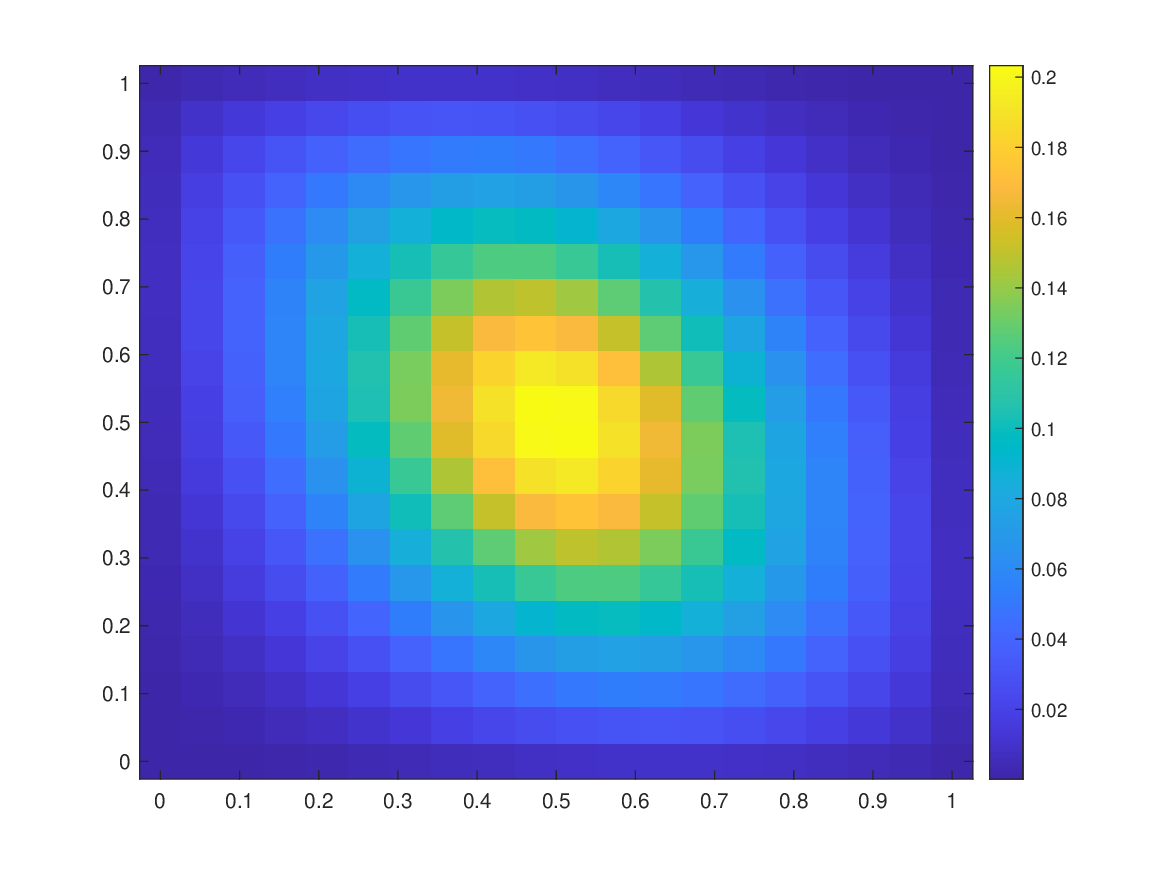}

\caption{Case 2. Top-Left: reference average solution in $\Omega_1$.  Top-Right: homogenized average solution in $\Omega_1$.  Bottom-Left: reference average solution in $\Omega_2$.  Bottom-Right: homogenized average solution in $\Omega_2$.}
\label{fig:case2_compare}
\end{figure}

For the second case, we change the permeability field to the one
shown in Figure \ref{fig:case2_perm}. 
We present $e_2^{(i)}$ in Table \ref{tab:case2}.
 Again, we observe that the proposed
approach provides an accurate approximation of the averaged solution
as we decrease the mesh size.
In Figure \ref{fig:case2_compare}, we depict upscaled solutions
and corresponding averaged fine-scale solutions. We observe a
 good agreement between coarse- and fine-grid solutions.
In the first table, we decrease the coarse-mesh size and the
period size at the same time. We observe that the error decreases
as the mesh size decreases.
Here, by choosing an appropriate number of layers, we observe that
the error remains small. 
In the second table,  we observe the convergence as 
we decrease the mesh size and fix
$\epsilon$. Again, the error decreases as we decrease the mesh size.

\begin{table}
\centering
\begin{tabular}{|c|c|c|c|}
\hline 
$H$ & $\epsilon$ & $e_{2}^{(1)}$ & $e_{2}^{(2)}$\tabularnewline
\hline 
$1/10$ & $1/10$ & $11.74\%$ & $1.61\%$\tabularnewline
\hline 
$1/20$ & $1/20$ & $4.18\%$ & $0.97\%$\tabularnewline
\hline 
$1/40$ & $1/40$ & $1.86\%$ & $1.08\%$\tabularnewline
\hline 
\end{tabular}
\begin{tabular}{|c|c|c|c|}
\hline 
$H$ & $\epsilon$ & $e_{2}^{(1)}$ & $e_{2}^{(2)}$\tabularnewline
\hline 
$1/10$ & $1/10$ & $11.74\%$ & $1.61\%$\tabularnewline
\hline 
$1/10$ & $1/20$ & $9.12\%$ & $6.13\%$\tabularnewline
\hline 
$1/10$ & $1/40$ & $8.05\%$ & $7.27\%$\tabularnewline
\hline 
\end{tabular}
\caption{Error comparison for Case 2.}
\label{tab:case2}
\end{table}

\section{Conclusions}

In this paper, we propose a general framework for 
multicontinuum homogenization. The method introduces several macroscopic
variables at each macroscale point using characteristic functions associated
with subdomains. The homogenization expansion is written
using macroscale variables and associated local cell problems. The local
cell problems are formulated as constraint problems in oversampled regions.
A decay of local cell solutions is needed for accurate approximations. This is not 
an easy task, in general, since the constraints are formulated in a spatially
localized fashion. We present an example of a mixed formulation of the elliptic equation, where we use some special formulations for cell problems. The proposed general framework shows that one can obtain various macroscale equations. We briefly discuss the relation to mixture theories.

\bibliographystyle{abbrv}
\bibliography{references,references4,references1,references2,references3,decSol}

\end{document}